\title{Fast winning strategies in Avoider-Enforcer games}
\author{{Dan Hefetz \thanks{Institute of Theoretical Computer Science, ETH Zurich,
CH-8092 Switzerland; and School of Computer Science,
Raymond and Beverly Sackler Faculty of Exact Sciences, Tel Aviv
University, Tel Aviv, 69978, Israel. Email: dan.hefetz@inf.ethz.ch.
This paper is a part of the author's Ph.D. under the supervision of
Prof.\ Michael Krivelevich.}} \quad {Michael Krivelevich
\thanks{School of Mathematical Sciences, Raymond and Beverly
Sackler Faculty of Exact Sciences, Tel Aviv University, Tel Aviv,
69978, Israel. Email: krivelev@post.tau.ac.il. Research supported in
part by a USA-Israel BSF grant, a grant from the Israel Science
Foundation and by a Pazy Memorial Award.}} \quad {Milo\v{s} Stojakovi\'{c}
\thanks{Department of Mathematics and Informatics, University
of Novi Sad, Serbia. Email: smilos@inf.ethz.ch. Partially supported
by the Ministry of Science, Republic of Serbia, and Provincial
Secretariat for Science, Province of Vojvodina.}} \quad {Tibor
Szab\'o
\thanks{Institute of Theoretical Computer Science, ETH Zurich,
CH-8092 Switzerland. Email: szabo@inf.ethz.ch.}}}

\documentclass[12pt]{article}
\usepackage{amsmath,amssymb,latexsym,color,epsfig}

\newif\ifnotesw\noteswtrue% T to show comments F supresses.

\def\cupdis{\dot \cup}
\def\ce{{\cal E}}
\def\ca{{\cal A}}

\def\dd{\overline{\overline{d}}}
\def\ddf{\lfloor\overline{\overline{d}}\rfloor}

\newtheorem{theorem}{Theorem}%[section]

\newtheorem{observation}[theorem]{Observation}

\newtheorem{question}[theorem]{Question}

\def\cf{{\mathcal F}}

\renewcommand{\epsilon}{\varepsilon}

\begin{document}
\maketitle

\begin{abstract}
In numerous positional games the identity of the winner is easily
determined. In this case one of the more interesting questions is not {\em who}
wins but rather {\em how fast} can one win. These type of problems were
studied earlier for Maker-Breaker games; here we initiate their study
for unbiased Avoider-Enforcer games played on the edge set
of the complete graph $K_n$ on $n$ vertices. For several games that
are known to be an Enforcer's win, we estimate quite precisely the minimum
number of moves Enforcer has to play in order to win.
We consider the non-planarity game, the connectivity game and the
non-bipartite game.
%It turns out that for these three games one can
%determine, up to a factor of two, the leading term of the minimum number of
%moves Enforcer needs in order to win, just by looking at extremal properties
%of the corresponding families of winning sets. Therefore, our main goal is to find the
%constant in the leading term of this value, as well as to estimate
%the lower order terms as accurately as possible.
\end{abstract}

\section{Introduction}
\label{sec::intro}

Let ${\mathcal F}$ be a hypergraph. In an unbiased Avoider-Enforcer
game ${\mathcal F}$ two players, called Avoider and Enforcer, take
turns selecting previously unclaimed vertices of $\mathcal{F}$,
with Avoider going first. Each player selects one vertex per turn,
until all vertices are claimed. Enforcer wins if Avoider claims all
the vertices of some hyperedge of $\mathcal{F}$; otherwise Avoider
wins. We refer to the family of hyperedges of $\cf$ as the family of losing
sets. In this paper our attention is restricted to games which are
played on the edges of the complete graph on $n$ vertices, that is,
the vertex set of ${\mathcal F}$ will always be $E(K_n)$.

Many positional games that were previously studied, are known to be
an easy win for Enforcer (for a comprehensive reference on
positional games the reader is referred to~\cite{Beck2}). For
example, the non-planarity game, where the goal of Avoider is to
avoid a non-planar graph, exhibits that kind of behavior -- Avoider
creates a non-planar graph and thus loses the game in the end,
irregardless of his strategy, the simple reason being that every
graph on $n$ vertices with more than $3n-6$ edges is non-planar.
Thus, for games of this type, a more interesting question to ask is
not {\em who wins} but rather \emph{how long does it take} the
winner to reach a winning position. This is the general problem we
address in this paper. To the best of our knowledge, ``fast
winning'' in Avoider-Enforcer games has not been studied before this
paper. On the other hand, there are quite a few results concerning
the analogous notion for Maker-Breaker games (see,
e.g.,~\cite{Beckfast,Bednarska,HKSS2,Pekec}).

For a hypergraph ${\mathcal F}$, let $\tau_E({\mathcal F})$ denote
the smallest integer $t$ such that Enforcer has a strategy to win
the game on ${\mathcal F}$ within $t$ moves. For the sake of
completeness, we define $\tau_E({\cal F}) = \infty$ if the game is
an Avoider's win.

One general way to approach the problem of determining the threshold
$\tau_E(\cf)$ is by investigating the extremal properties of the hypergraph
$\cf$. For convenience, let us assume that the set of hyperedges of
$\cf$ is a monotone increasing family of sets. If this is not the
case, we can extend it to an increasing family by adding all the
supersets of its elements -- this operation clearly does not change
the outcome of the game. The {\em extremal number} (or Tur\'an number)
of the hypergraph $\cf$ is defined by
$$
\text{ex}(\cf) = \max\left\{|A|: A\subseteq V(\cf),\,\,A\not\in
E(\cf) \right\}.
$$
Then the minimum move number $\tau_E (\cf )$ can be determined up to a
factor of two.
\begin{observation} \label{obs}
Given a monotone increasing family ${\mathcal F}$ of hyperedges, we have
$$
\frac{1}{2}\text{\rm ex}(\cf)+1 \leq  \tau_E(\cf) \leq
\text{\rm ex}(\cf) +1.
$$
\end{observation}
{\em Proof.}
To prove the lower bound, let Avoider fix an arbitrary subset
$A$ of $V(\cf)$ before the game starts, such that $A$ is not an edge of
$\cf$ and $|A|=\text{ex}(\cf )$. Then, during the game, Avoider
just claims elements of $A$ for as long as possible.
This way he will be able to claim at least half of the elements of
$A$ without losing.

For the upper bound, observe that Enforcer will surely win after
$\text{ex}(\cf) +1$ rounds irregardless of his strategy. Indeed, at
that point, Avoider has claimed $\text{ex}(\cf) +1$ vertices, and
every set with that many vertices is an edge of $\cf $. \hfill
$\Box$

\subsection{Our results}

As we have already mentioned, in the Avoider-Enforcer non-planarity
game Avoider loses the game as soon as his graph becomes non-planar.
The biased version of this game was studied in~\cite{HKSS}. Denote
by ${\cal NP}_n$ the hypergraph whose hyperedges are the edge-sets
of all non-planar graphs on $n$ vertices. From
Observation~\ref{obs}, we obtain
$$
\frac{3}{2}n-2\leq\tau_E({\cal NP}_n) \leq 3n-5.
$$
The following theorem asserts that this upper bound is essentially
tight, that is, Avoider can refrain from building a non-planar graph
for at least $(3-o(1))n$ moves. More precisely,
\begin{theorem}
\label{planarity}
$$\tau_E({\cal NP}_n) > 3n - 28\sqrt{n}.$$
\end{theorem}

In the Avoider-Enforcer non-bipartite game Avoider loses the game as
soon as his graph first becomes non-bipartite. Clearly, this game is
equivalent to the game in which Avoider's goal is to avoid creating
an odd cycle. Denote by ${\mathcal NC}_n^2$ the hypergraph whose
hyperedges are the edge-sets of all non-bipartite graphs on $n$
vertices. Mantel's Theorem asserts that the bipartite graph on $n$
vertices which maximizes the number of edges is the complete
bipartite graph with a balanced partition. Hence, it follows from
Observation~\ref{obs} that
$$
\frac{1}{2}\left\lfloor\frac{n^2}{4}\right\rfloor +1\leq \tau_E({\cal
NC}_n^2)\leq \left\lfloor\frac{n^2}{4}\right\rfloor +1.
$$
In the next theorem we improve the upper bound substantially and establish
that the lower bound is asymptotically correct. We also slightly improve
the lower bound and thus determine
the order of magnitude of the second order term of $\tau_E({\cal
NC}_n^2)$.

\begin{theorem}
\label{coloring} $$\tau_E({\cal NC}_n^2) = \frac{n^2}{8} + \Theta(n).$$
\end{theorem}

%\begin{theorem}
%\label{coloring} $$\frac{n^2}{8} + \frac{n-2}{12} \leq \tau_E({\cal
%NC}_n^2) \leq \frac{n^2}{8} + \frac{n}{2} + 1.$$
%\end{theorem}

Note that the non-bipartite game is just a special case of the
non-$k$-colorability game ${\mathcal NC}_n^k$, where Avoider loses the game
as soon as his graph becomes non-$k$-colorable.
Observation~\ref{obs} can be readily applied, but it would be
interesting to obtain tighter bounds, as in the case $k=2$.

%Avoider can play for at least $(1-o(1))\frac{(k-1)n^2}{4k}$ moves
%without losing by simply fixing a copy of the $k$-partite
%Tur\'an-graph and claiming half of its edges. On the other hand, it
%is not hard to see that the game is an Enforcer's win if it is
%played until the end (see~\cite{HKSS}), so Avoider will lose after
%at most $\frac{1}{2}{n \choose 2} \approx \frac{n^2}{4}$ moves.

Finally, we consider two Avoider-Enforcer games that turn out to be of
similar behavior. In the positive min-degree game, Enforcer wins as soon as
the minimum degree in Avoider's graph becomes positive, and in the
connectivity game, Enforcer wins as soon as Avoider's graph becomes
connected and spanning. Denote by ${\mathcal D}_n$ and ${\mathcal T}_n$
the hypergraphs whose hyperedges are the edge-sets of all graphs with a
positive minimum degree, and the edge-sets of all graphs that are
connected and spanning, respectively.

Clearly, we have $\tau_E({\mathcal D}_n) \leq \tau_E({\mathcal T}_n)$,
since ${\mathcal D}_n \supseteq {\mathcal T}_n$. As
$\text{\rm ex}({\cal D}_n) = \text{\rm ex}({\cal T}_n) = {n-1
\choose 2}$, Observation~\ref{obs} implies
$$
\frac{1}{2}{n-1 \choose 2} +1\leq \tau_E({\cal D}_n) \leq
\tau_E({\cal T}_n) \leq {n-1 \choose 2} +1.
$$
Moreover, as Enforcer wins both games (see~\cite{HKS}), we have
$$\tau_E({\cal D}_n), \tau_E({\cal T}_n) \leq \frac{1}{2}{n \choose
2},
$$
which determines both parameters asymptotically and shows that they
are ``quite close to each other''. This is somewhat
reminiscent of the well-known property of random graphs, that the
hitting time of being connected and the hitting time of having
minimum positive degree are a.s.\ the same, and it motivates us to
raise the following question.
\begin{question} \label{hittingtime} Is it true that
$\tau_E({\cal D}_n) = \tau_E({\cal T}_n)$ holds for sufficiently large $n$?
\end{question}

The following theorem can be considered as a first step
towards an affirmative answer to Question~\ref{hittingtime}. We improve the aforementioned  lower and upper bounds, determining in the process the second order term and the
order of magnitude of the third for both of these parameters.
\begin{theorem} \label{posdeg=connect}
$$
\tau_E({\cal D}_n), \tau_E({\cal T}_n) = \frac{1}{2}{n-1
\choose 2}+ \Theta\left( \log n\right).
$$
\end{theorem}

%\begin{theorem} \label{posdeg=connect}
%$$
%\frac{1}{2}{n-1 \choose 2} + \left(\frac{1}{4} - o(1)\right)\log n <
%\tau_E({\cal D}_n) \leq \tau_E({\cal T}_n) \leq \frac{1}{2}{n-1
%\choose 2}+ 2\log_2 n + 1.
%$$
%\end{theorem}

%\noindent The rest of the paper is organized as follows: In
%Section~\ref{sec::avoiderenforcer} we prove
%Theorems~\ref{planarity},~\ref{coloring} and~\ref{posdeg=connect},
%and in Section~\ref{sec::openprobs} we present some open
%problems.

%\subsection{Preliminaries}

\vspace{1cm}

For the sake of simplicity and clarity of presentation, we
omit floor and ceiling signs whenever these are not crucial. Some of
our results are asymptotic in nature and, whenever necessary, we
assume that $n$ is sufficiently large. Throughout the paper, $\log$
stands for the natural logarithm. Our graph-theoretic notation is
standard and follows that of~\cite{Diestel}.

%In particular, we use the following: for a graph $G$, denote its set
%of vertices by $V(G)$, and its set of edges by $E(G)$. Moreover, let
%$v(G) = |V(G)|$ and $e(G) = |E(G)|$. For a graph $G=(V,E)$ and a set
%$A \subseteq V$ denote by $G[A]$ the subgraph of $G$ induced by $A$.
%Let $N_G(A) = \{u \in V : \exists w \in A, (u,w) \in E\}$ be the
%neighborhood of $A$ in $G$ and let $\Gamma_G(A) = N_G(A) \setminus
%A$ be the external-neighborhood of $A$ in $G$. Sometimes, when there
%is no risk of confusion, we abbreviate $N_G(A)$ to $N(A)$ and
%$\Gamma_G(A)$ to $\Gamma(A)$. For a hypergraph $\cf$, we denote its
%set of vertices by $V(\cf)$, and its set of hyperedges by $E(\cf)$.

\section{The strategies} \label{sec::avoiderenforcer}

\subsection{Keeping the graph planar for long}\label{sec::planarity}
{\em Proof of Theorem~\ref{planarity}}
We begin by introducing some terminology. Let
$v$ be a vertex, and let $S$ be a set of vertices. Let $N_\ca(v,S)$
denote the set of neighbors of $v$ in Avoider's graph, belonging to
$S$. Similarly, let $N_\ce(v,S)$ denote the set of neighbors of $v$
in Enforcer's graph, belonging to $S$.

We will provide Avoider with a strategy for keeping his graph
planar for at least $3n-28\sqrt{n}$ rounds. The strategy consists of
three stages.

Before the game starts, we partition the vertex set
$$
V(K_n) = \{v_1\}\cupdis \{v_2\} \cupdis A \cupdis N_{1,1} \cupdis
N_{1,2}\cupdis N_{2,1} \cupdis N_{2,2},
$$
such that $|N_{1,1}|=|N_{1,2}|=|N_{2,1}|=|N_{2,2}|=\sqrt{n}-1$ and
$|A|=n- 4\sqrt{n}+2$.

In the first stage, Avoider claims edges according to a simple
pairing strategy. For every vertex $a \in A$, we pair up the edges
$(a,v_1)$ and $(a,v_2)$. Whenever Enforcer claims one of the paired
edges, Avoider immediately claims the other edge of that pair. If
Enforcer claims an edge which does not belong to any pair, then
Avoider claims the edge $(a,v_1)$, for some $a \in A$, for which
neither $(a,v_1)$ nor $(a,v_2)$ were previously claimed. He then
removes the pair $(a,v_1),(a,v_2)$ from the set of considered edge
pairs.

The first stage ends as soon as Avoider connects every $a \in A$ to
either $v_1$ or $v_2$. Note that, at that point, Avoider's graph
consists of two vertex-disjoint stars centered at $v_1$ and $v_2$,
and the isolated vertices in $N_{1,1} \cup N_{1,2} \cup N_{2,1} \cup
N_{2,2}$. Hence, during the first stage, Avoider has claimed exactly
$n- 4\sqrt{n}+2$ edges. Define $A_1:=N_\ca (v_1,A)$, and $A_2:=N_\ca
(v_2,A)$.

Before the second stage starts, we pick four vertices $n_{1,1}\in
N_{1,1}$, $n_{1,2}\in N_{1,2}$, $n_{2,1}\in N_{2,1}$ and $n_{2,2}\in
N_{2,2}$, such that $|N_\ce (n_{i,j}, A)|\leq \sqrt{n}$, for every
$i,j\in\{1,2\}$. Clearly, such a choice of vertices is possible as
the total number of edges Enforcer has claimed during the first
stage is $n-4\sqrt{n}+2< \sqrt{n}\cdot (\sqrt{n}-1)$. Define
$G_1:=N_\ce (n_{1,1},A_1)\cup N_\ce (n_{1,2},A_1)$, and $G_2:=N_\ce
(n_{2,1},A_2)\cup N_\ce (n_{2,2},A_2)$. Note that $|G_1|\leq
2\sqrt{n}$, $|G_2|\leq 2\sqrt{n}$, and $|N_\ce (n_{1,1},
A_1\setminus G_1)| = |N_\ce (n_{1,2}, A_1\setminus G_1)| = |N_\ce
(n_{2,1}, A_2\setminus G_2)| = |N_\ce (n_{2,2}, A_2\setminus
G_2)|=0$.

Using a pairing strategy similar to the one used in the first
stage, Avoider connects each vertex of $A_1\setminus G_1$ to either
$n_{1,1}$ or $n_{1,2}$, and each vertex of $A_2\setminus G_2$ to
either $n_{2,1}$ or $n_{2,2}$. More precisely, for every $a\in
A_1\setminus G_1$ we pair up the edges $(a,n_{1,1})$ and
$(a,n_{1,2})$, and for every $a\in A_2\setminus G_2$ we pair up
edges $(a,n_{2,1})$ and $(a,n_{2,2})$. Avoider then proceeds as in
the first stage.

The second stage ends as soon as Avoider connects every $a \in A_1
\setminus G_1$ to either $n_{1,1}$ or $n_{1,2}$, and every $a \in
A_2 \setminus G_2$ to either $n_{2,1}$ or $n_{2,2}$. We define
$A_{1,1}:= N_\ca (n_{1,1},A_1)$, $A_{1,2}:= N_\ca (n_{1,2},A_1)$,
$A_{2,1}:= N_\ca (n_{2,1},A_2)$ and $A_{2,2}:= N_\ca (n_{2,2},A_2)$.
Since $|A_{1,1}|+|A_{1,2}| = |A_1|-|G_1|$, $|A_{2,1}|+|A_{2,2}| =
|A_2|-|G_2|$ and $|A_1|+|A_2| = |A|$, we infer that the number of
edges Avoider has claimed in the second stage is at least
$n-8\sqrt{n}$. Note that during the first two stages Avoider did
not claim any edge with both endpoints in one of the sets
$A_{1,1}$, $A_{1,2}$, $A_{2,1}$, $A_{2,2}$.

In the third stage, Avoider claims only edges with both endpoints
contained in the sets $A_{i,j}$, for some $i,j\in\{1,2\}$. His goal
in this stage is to build a ``large'' linear forest in $A_{1,1}$. (A
\emph{linear forest} is a vertex-disjoint union of paths.) In the
beginning of the third stage, Avoider's graph induced on the
vertices of $A_{1,1}$ is empty, that is, it consists of $|A_{1,1}|$
paths of length 0 each. For as long as possible, Avoider claims
edges that connect endpoints of two of his paths in $A_{1,1}$,
creating a longer path. When this is no longer possible, every edge
that connects endpoints of two different paths must have been
previously claimed by Enforcer. Since the total number of edges that
Enforcer has claimed so far is at most $3n$, the number of paths of
Avoider in $A_{1,1}$ is at most $2\sqrt{n}$. Hence, Avoider has
claimed at least $|A_{1,1}| - 2\sqrt{n}$ edges to this point of the
third stage.
\begin{figure}[tbp]
  \centering %
  \epsfig{file=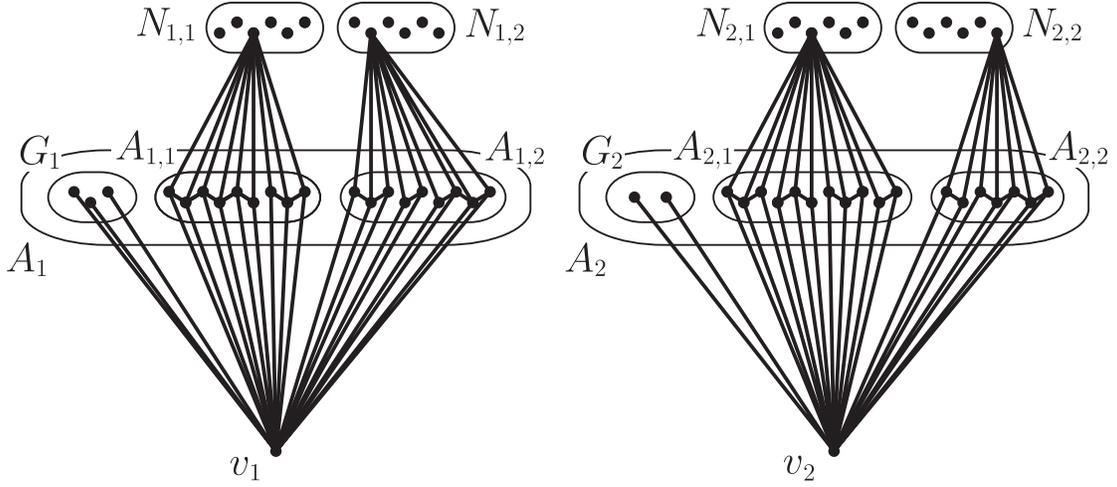, scale=1.4}
  \caption{Avoider's graph.\label{fig1}}
\end{figure}

Similarly, Avoider builds a ``large'' linear forest in $A_{1,2}$,
$A_{2,1}$, and finally $A_{2,2}$, all in the same way. Thus, the
total number of edges he claims during the third stage is at
least
\begin{eqnarray*}
\sum_{i,j\in\{1,2\}} (|A_{i,j}|-2\sqrt{n})
&\geq& |A_1|-|G_1| + |A_2|-|G_2| - 8 \sqrt{n} \\
&\geq& |A| - 12\sqrt{n} \\
&\geq& n - 16 \sqrt{n}.
\end{eqnarray*}

The total number of edges claimed by Avoider during the entire game
is therefore at least $(n - 4\sqrt{n}) + (n-8\sqrt{n}) + (n - 16
\sqrt{n}) = 3n - 28\sqrt{n}$. Moreover, at the end of the third
stage (which is also the end of the game), Avoider's graph is the
pairwise edge disjoint union of two stars, four other graphs - each
being a subgraph of a union of $K_{2,n_i}$ and a linear forest which
is restricted to one side of the bipartition (see
Figure~\ref{fig1}). Clearly, such a graph is planar. {\hfill $\Box$}

\medskip

\subsection{Forcing and avoiding odd cycles}
\label{sec::coloring}

{\em Proof of Theorem~\ref{coloring}}

\paragraph{Forcing an odd cycle fast.}
First, we provide Enforcer with a strategy that will force Avoider
to claim the edges of an odd cycle during the first $\frac{n^2}{8} +
\frac{n}{2} + 1$ moves. In every stage of the game, each connected
component of Avoider's graph is a bipartite graph with a unique
bipartition of the vertices (we stop the game as soon as Avoider is
forced to close an odd cycle). In every move, Enforcer's primary
goal is to claim an edge which connects two opposite sides of the
bipartition of one of the connected components of Avoider's graph.
If no such edge is available, then Enforcer claims an arbitrary
edge, and that edge is marked as ``possibly bad''. Clearly, in his
following move Avoider cannot play inside any of the connected
components of his graph either, and so he is forced to merge two of
his connected components (that is, he has to claim an edge $(x,y)$
such that $x$ and $y$ are in different connected components of his
graph). As the game starts with $n$ connected components, this
situation can occur at most $n-1$ times.

Therefore, when Avoider is not able to claim any edge without
creating an odd cycle, his graph is bipartite, and all of Enforcer's
edges, except some of the ``possibly bad'' ones, are compatible with
the bipartition of Avoider's graph. The total number of edges that
were claimed by both players to this point is at most $\frac{n^2}{4}
+ n - 1$, and so the total number of moves Avoider has played in the
entire game is at most $\frac{n^2}{8} + \frac{n}{2} + 1$.

\paragraph{Avoiding odd cycles for long.}
Next, we provide Avoider with a strategy for keeping his graph
bipartite for at least $\frac{n^2}{8} + \frac{n-2}{12}$ rounds. For
technical reasons we assume that $n$ is even; however, a similar
statement holds for odd $n$ as well. During the game Avoider will
maintain a family of ordered pairs $(V_1, V_2)$, where $V_1, V_2
\subseteq V(K_n)$, $V_1 \cap V_2 = \emptyset$ and $|V_1|=|V_2|$,
which he calls \emph{bi-bunch}es. We say that two bi-bunches $(V_1,
V_2)$ and $(V_3, V_4)$ are disjoint if $(V_1 \cup V_2) \cap (V_3
\cup V_4) = \emptyset$. At any point of the game, Avoider calls a
vertex \emph{untouched} if it does not belong to any bi-bunch and
all the edges incident with it are unclaimed. During the entire
game, we will maintain a partition of the vertex set $V(K_n)$ into a
number of pairwise disjoint bi-bunches, and a set of untouched
vertices.

Avoider starts the game with $n$ untouched vertices and no
bi-bunches. In every move, his primary goal is to claim an edge
\emph{across} some existing bi-bunch, that is, an edge $(x,y)$ where
$x \in V_1$ and $y \in V_2$ for some bi-bunch $(V_1,V_2)$. If no
such edge is available, then he claims an edge joining two untouched
vertices $x$ and $y$, introducing a new bi-bunch $(\{x\},\{y\})$. If
he is unable to do that either, then he claims an edge connecting
two bi-bunches, that is, an edge $(x,y)$ such that there exist two
bi-bunches $(V_1, V_2)$ and $(V_3, V_4)$ with $x \in V_1$ and $y \in
V_3$. He then replaces these two bi-bunches with a single new one
$(V_1 \cup V_4, V_2 \cup V_3)$.

Whenever Enforcer claims an edge $(x,y)$ such that neither $x$ nor
$y$ belong to any bi-bunch, we introduce a new bi-bunch $(\{x,y\},
\{u,v\})$, where $u$ and $v$ are arbitrary untouched vertices. If at
that point of the game there are no untouched vertices (clearly this
can happen at most once), then the new bi-bunch is just
$(\{x\},\{y\})$.
%Danny: WE HAVE TO TAKE THIS CASE INTO
%ACCOUNT IN OUR ANALYSIS.
%Milos: I think that this is included if we assume that n/6 is
%rounded down.
If Enforcer claims an edge $(x,y)$ such
that there is a bi-bunch $(V_1,V_2)$ with $x \in V_1$ and $y$ is
untouched, then the bi-bunch $(V_1,V_2)$ is replaced with $(V_1 \cup
\{y\}, V_2 \cup \{u\})$, where $u$ is an arbitrary untouched vertex.
Finally, if Enforcer claims an edge $(x,y)$ such that there are
bi-bunches $(V_1,V_2)$ and $(V_3,V_4)$ with $x \in V_1$ and $y \in
V_3$, than these two bi-bunches are replaced with a single one $(V_1
\cup V_3, V_2 \cup V_4)$. Note that by following his strategy, and
updating the bi-bunch partition as described, Avoider's graph will
not contain an edge with both endpoints in the same side of a
bi-bunch at any point of the game.

Observe that the afore-mentioned bi-bunch maintenance rules imply
the following. If Enforcer claims an edge $(x,y)$, such that before
that move $x$ was an untouched vertex, then the edge $(x,y)$ will be
contained in the same side of some bi-bunch, that is, after that
move there will be a bi-bunch $(V_1,V_2)$ with $x,y \in V_1$ (unless
$x$ and $y$ were the last two isolated vertices).

Assume that in some move Avoider claims an edge $(x,y)$, such that
before that move $x$ was an untouched vertex. It follows from
Avoider's strategy that $y$ was untouched as well, and there were no
unclaimed edges across a bi-bunch at that point. Thus, in his next
move, Enforcer will also be unable to claim an edge across a
bi-bunch and so, by the bi-bunch maintenance rules for Enforcer's
moves, the edge he will claim in that move will have both its
endpoints in the same side of some bi-bunch.

By the previous paragraphs, we conclude that after every round in
which at least one of the players claims an edge which is incident
with an untouched vertex (which is not the next to last untouched
vertex), the edge Enforcer claims in this round will be contained in
the same side of some bi-bunch. By the bi-bunch maintenance rules,
during every round the number of untouched vertices is decreased by
at most 6. Hence, by the time all but two vertices are not untouched
at least $(n-2)/6$ edges of Enforcer will be contained in the same
side of a bi-bunch. Therefore, when Avoider can no longer claim an
edge without creating an odd cycle, both players have claimed
together all the edges of a balanced bipartite graph which is in
compliance with the bi-bunch bipartition, and at least another
$(n-2)/6$ edges. This gives a total of at least $\frac{n}{2} \cdot
\frac{n}{2} + (n-2)/6$ edges claimed, which means that at least
$\frac{n^2}{8} + \frac{n-2}{12}$ rounds were played to that point.
{\hfill $\Box$}

\medskip

\subsection{Spanning trees and isolated vertices}
\label{sec::posdeg=connect}

{\em Proof of Theorem~\ref{posdeg=connect}.} Clearly $\tau_E({\cal
D}_n) \leq \tau_E({\cal T}_n)$ and so it suffices to prove that
$\tau_E({\cal T}_n) \leq \frac{1}{2}{n-1 \choose 2} + 2\log_2 n + 1$
and that, $\tau_E({\cal D}_n) > \frac{1}{2}{n-1 \choose 2} + (1/4 -
\varepsilon) \log n$, for every $\varepsilon>0$ and sufficiently
large $n$.

\paragraph{Forcing a spanning tree fast.}
Starting with the former inequality, we provide Enforcer with a
strategy to force Avoider to build a connected spanning graph within
$\frac{1}{2}{n-1 \choose 2}+ 2\log_2 n + 1$ rounds. At any point of
the game, we call an edge that was not claimed by Avoider
\emph{safe}, if both its endpoints belong to the same connected
component of Avoider's graph. An edge which is not safe and was not
claimed by Avoider is called \emph{dangerous}. Denote by $G_D$ the
graph consisting of dangerous edges claimed by Enforcer. We will
provide Enforcer with a strategy to make sure that, throughout the
game, the maximum degree of the graph $G_D$ does not exceed $4k$,
where $k=\log_2 n$.

Assuming the existence of such a strategy, the assertion of the
theorem readily follows. Indeed, assume for the sake of
contradiction that after $\frac{1}{2}{n-1 \choose 2}+ 2\log_2 n + 1$
rounds have been played (where Enforcer follows the afore-mentioned
strategy), Avoider's graph is disconnected. Let $C_1, \ldots, C_r$,
where $r \geq 2$ and $|C_1| \leq \ldots \leq |C_r|$, be the
connected components in Avoider's graph at that point. By Enforcer's
strategy, the maximum degree of the graph $G_D$ does not exceed
$4k$. Hence, the number of edges claimed by both players to this
point does not exceed
$$
\sum_{i=1}^{r}{|C_i| \choose 2} + 4k \sum_{i=1}^{r-1}|C_i|.
$$
Assuming that $r \geq 2$ and $n$ is sufficiently large, this sum
above attains its maximum for $r=2$, $|C_1| = 1$ and $|C_2| = n-1$;
that is, the sum is bounded from above by ${n-1 \choose 2}+ 4 \log_2
n$ - a contradiction.

Now we provide Enforcer with a strategy for making sure that,
throughout the game, the maximum degree of the graph $G_D$ does not
exceed $4k$. In every move, if there exists an unclaimed safe edge,
Enforcer claims it (if there are several such edges, Enforcer claims
one arbitrarily). Hence, whenever Enforcer claims a dangerous edge,
Avoider has to merge two connected components of his graph in the
following move, and the number of Avoider's connected components is
decreased by one. We will use this fact to estimate the number of
dangerous edges at different points of the game.

When all edges within each of the connected components of Avoider's
graph are claimed, Enforcer has to claim a dangerous edge. His
strategy for choosing dangerous edges is divided into two phases.
The first phase is divided into $k$ stages. In the $i$th stage
Enforcer will make sure that the maximum degree of the graph $G_D$
is at most $2i$; other than that, he claims dangerous edges
arbitrarily. He proceeds to the following stage only when it is not
possible to play in compliance with this condition. Let $c_i$, $i=1,
\dots, k$, denote the number of connected components in Avoider's
graph after the $i$th stage. Let $c_0 = n$, be the number of
components at the beginning of the first stage. During the $i$th
stage, a vertex $v$ is called \emph{saturated}, if $d_{G_D}(v)=2i$.
Note that at the beginning of the first stage the maximum degree of
$G_D$ is $2\cdot 0 = 0$.

We will prove by induction that $c_i \leq n 2^{-i} + 2i$, for all
$i=0, 1,\dots,k$. The statement trivially holds for $i=0$.

% In the beginning of the first stage Avoider's graph has $n$
% connected components, and in the end of the first stage Avoider's
% graph has $c_1$ connected components. It follows that the number of
% dangerous edges that Enforcer has claimed during the first stage
% is at most $n - c_1$. Hence, the total number of saturated vertices
% cannot be greater than $n-c_1$. Clearly, two vertices that are not
% saturated cannot be in different connected components of Avoider's
% graph, as otherwise Enforcer would have claimed the edge that
% connects them, contrary to the fact that the first stage has ended.
% It follows that there is a component of size at least $c_1$. As the
% number of components is also $c_1$, we get $(c_1 - 1) + c_1 \leq
% n$ entailing $c_1 \leq (n+1)/2 \leq n/2 + 1$.

Next, assume that $c_{j} \leq n 2^{-j} + 2j$, for some $0 \leq j <
k$. At the beginning of the $(j+1)$st stage Avoider's graph has
exactly $c_j$ connected components, and at the end of this stage it
has exactly $c_{j+1}$ components. It follows that during this stage
Avoider merged two components of his graph $c_j - c_{j+1}$ times.
Hence, Enforcer has not claimed more than $c_{j}-c_{j+1}$ dangerous
edges during the $(j+1)$st stage. As the maximum degree of the graph
$G_D$ before this stage was $2j$, the number of saturated vertices
at the end of the $(j+1)$st stage is at most $c_{j}-c_{j+1}$. It
follows that there are at least $n - (c_j-c_{j+1})$ non-saturated
vertices at this point. The non-saturated vertices must be covered
by at most $2(j+1)$ connected components of Avoider's graph. Indeed,
assume for the sake of contradiction that there are non-saturated
vertices $u_1,u_2, \ldots, u_{2j+3}$ and connected components
$U_1,U_2, \ldots, U_{2j+3}$, such that $u_p \in U_p$ for every $1
\leq p \leq 2j+3$. Since $deg_{G_D}(u_p) \leq 2j+1$ for every $1
\leq p \leq 2j+3$, it follows that there must exist an unclaimed
edge $(u_r,u_s)$ for some $1 \leq r < s \leq 2j+3$, contradicting
the fact that the $(j+1)$st stage is over. Therefore, there are at
least $c_{j+1} - 2(j+1)$ connected components in Avoider's graph
that do not contain any non-saturated vertex. Clearly every such
component has size at least one, entailing $(c_{j+1} - 2j - 2) +
(n-c_{j} + c_{j+1}) \leq n$. Applying the inductive hypothesis we
get $c_{j+1} \leq c_j/2 + j+1 \leq n 2^{-(j+1)} + 2(j+1)$. This
completes the induction step.

It follows, that at the end of the first phase, after the $k$th
stage, the number of connected components in Avoider's graph, is at
most $c_k \leq n 2^{-k} + 2k \leq 2k+1$.

In the second phase, whenever Enforcer is forced to claim a
dangerous edge, he claims one arbitrarily. Since at the beginning of
the second phase, there are at most $2k+1$ connected components in
Avoider's graph, Enforcer will claim at most $2k$ dangerous edges
during this phase.

It follows that at the end of the game, the maximum degree in $G_D$
will be at most $4k$, as claimed.

\paragraph{Keeping an isolated vertex for long.}
Fix $\varepsilon > 0$ and set $l:=\frac{1-4\varepsilon}{2} \log n$.
We provide Avoider with a strategy to keep an isolated vertex in his
graph for at least $\frac{1}{2}{n-1 \choose 2} + \frac{l}{2}$
rounds.

Throughout the game, Avoider's graph will consist of one connected
component, which we denote by $C$, and $n - |C|$ isolated vertices.
A vertex $v \in V(K_n) \setminus C$ is called \emph{bad}, if there
is an even number of unclaimed edges between $v$ and $C$; otherwise,
$v$ is called \emph{good}.

For every vertex $v \in V(K_n)$ let $d_{\cal E}(v)$ denote the
degree of $v$ in Enforcer's graph. If at any point of the game there
exists a vertex $v \in V(K_n) \setminus C$ such that $d_{\cal E}(v)
\geq l$, then Avoider simply proceeds by arbitrarily claiming edges
which are not incident with $v$, for as long as possible. The total
number of rounds that will be played in that case is at least
$\frac{1}{2}{n-1 \choose 2}+ \frac{l}{2}$, which proves the theorem.
We will show that Avoider can make sure that such a vertex $v \in
V(K_n) \setminus C$, with $d_{\cal E}(v) \geq l$, will appear before
the order of his component $C$ reaches $n - l \varepsilon^{-1} - 1$.
Hence, from now on, we assume that $|C| \leq n - l \varepsilon^{-1}
- 2$.

Whenever possible, Avoider will claim an edge with both endpoints in
$C$. If this is not possible, he will join a new vertex to the
component, that is, he will connect it by an edge to an arbitrary
vertex of $C$. Note that this is always possible. Indeed, assume
that every edge between $C$ and $V(K_n) \setminus C$ was already
claimed by Enforcer. If $|C| \geq l$ then there exists a vertex $v
\in V(K_n)$ such that $d_{\cal E}(v) \geq l$ and so we are done by
the previous paragraph. Otherwise, $|C| < l$ and thus, until this
point, Enforcer has claimed at most $l^2 < l(n-l)$ edges. As for the
way he chooses this new vertex, we consider three cases. Let $\dd$
denote the average degree in Enforcer's graph, taken over all the
vertices of $V(K_n) \setminus C$, that is,
$$
\dd:=\frac{\sum_{v \in V(K_n) \setminus C} d_{\cal E}(v)}{n-|C|}.
$$
Throughout the case analysis, $C$ and $\dd$ represent the values as
they are just before Avoider makes his selection.

\begin{enumerate}
\item There exists a vertex $v \in V(K_n) \setminus C$, such that $d_{\cal E}(v) \leq
\dd - 1$.

Avoider joins $v$ to his component $C$. %Let $\ddo$ and $\ddn$ denote
%the values of $\dd$ before and after this move.
Then $|C|$ increases by one, and the new value of $\dd$ is at least
\begin{eqnarray*}
\frac{(n-|C|)\dd - (\dd -1)}{n-|C|-1} &=& \dd + \frac{1}{n-|C|-1}.
\end{eqnarray*}

\item Every vertex $v \in V(K_n) \setminus C$ satisfies $d_{\cal E}(v) >
\dd-1$, and $\dd < \ddf +1 - \varepsilon$.

Let $D$ denote the set of vertices $u \in V(K_n) \setminus C$ such
that $d_{\cal E}(u) = \ddf$. Note that there must be at least
$\varepsilon (n-|C|)$ vertices in $D$. We distinguish between the
following two subcases.
  \begin{enumerate}
  \item There is a good vertex in $D$. Avoider joins it to his
  component $C$ (if there are several good vertices, then he picks one
  arbitrarily).
  Since $v$ was a good vertex, Enforcer must claim at
  least one edge $(x,y)$ such that $x \notin C\cup\{ v\}$,
  before Avoider is forced again to join another vertex to his component.
  After this move of Enforcer $|C|$ is (still) increased by (just) one, and
  the new value of $\dd$ is at least
  \begin{eqnarray*}
  %&=& \frac{(n-|C|)\ddo - \ddf + d_{\cal E}(x)}{n-|C|-1}\\
  \frac{(n-|C|)\dd - \ddf + 1}{n-|C|-1}
  &\geq& \dd + \frac{1}{n-|C|-1}.
  \end{eqnarray*}

  \item All vertices in $D$ are bad.
  Knowing that $d_{\cal E}(v) \leq l-1$ for all vertices $v \in
  V(K_n)\setminus C$, and $|C| \leq n - l \varepsilon^{-1} - 2$, we have
  \begin{eqnarray*}
  \max_{v \in D} d_{\cal E}(v) = \ddf < l - 1 + 2\varepsilon
  \leq \varepsilon(n - |C|) - 1 \leq |D| - 1
  \end{eqnarray*}
  and hence there have to be two vertices $u,w \in D$ such that
  $(u,w)$ is unclaimed. Avoider
  joins $u$ to his component $C$, and thus $w$ becomes good.
  If Enforcer, in his next move, claims an edge $(w,v)$ for some $v
  \in C$,  then $|C|$ is increased by one and the
  new value of $\dd$ is at least
  \begin{eqnarray*}
  \frac{(n-|C|)\dd - \ddf + 1}{n-|C|-1}
  &\geq& \dd + \frac{1}{n-|C|-1}.
  \end{eqnarray*}

  Otherwise, in his next move Avoider joins $w$ to $C$. Since $w$
  was good, then, as in the previous subcase, Enforcer will be
  forced to claim an edge $(x,y)$ such that $x \notin C\cup \{ w\}$. After
  that move of Enforcer, we will have that $|C|$ is still increased
  just by two and the new value of $\dd$ is at least

  \begin{eqnarray*}
  %&=& \frac{(n-|C|)\ddo - \ddf - \ddf + d_{\cal E}(x)}{n-|C|-2}\\
  \frac{(n-|C|)\dd - \ddf - \ddf + 1}{n-|C|-2}
  &\geq& \dd + \frac{1}{n-|C|-2}.
  \end{eqnarray*}
  \end{enumerate}

\item Every vertex $v \in V(K_n) \setminus C$ satisfies $d_{\cal E}(v) >
\dd-1$, and $\dd \geq \ddf + 1 - \varepsilon$.

Let $D$ denote the set of vertices in $V(K_n) \setminus C$ with
degree either $\ddf$ or $\ddf +1$. Clearly, $|D| \geq
\frac{1}{2}(n-|C|)$. We distinguish between the following two
subcases.

  \begin{enumerate}
  \item There is a good vertex in $D$. Similarly to subcase 2(a),
  Avoider joins that vertex to his component $C$, and after
  Enforcer claims some edge with at least one endpoint outside $C$,
  we have that $|C|$ is increased by one and the new value of $\dd$ is at least
  \begin{eqnarray*}
  \frac{(n-|C|)\dd - (\dd + \varepsilon) + 1}{n-|C|-1}
  &=& \dd + \frac{1-\varepsilon}{n-|C|-1}.
  \end{eqnarray*}

  \item All vertices in $D$ are bad. Similarly to subcase 2(b),
  Avoider
  can find two vertices in $D$ such that the edge between them is
  unclaimed.
  He joins them to his component $C$, one after the other. After
  Enforcer claims some edge with at least one endpoint outside $C$,
  we have that $|C|$ increased by two and the new value of $\dd$ is at least
  \begin{eqnarray*}
  \frac{(n-|C|)\dd - (\dd + \varepsilon) - (\dd + \varepsilon) +
  1}{n-|C|-2}
  &=& \dd + \frac{1-2\varepsilon}{n-|C|-2}.
  \end{eqnarray*}
  \end{enumerate}

\end{enumerate}

It follows that in all cases the value of $\dd$ grows by at least
$\frac{1-2\varepsilon}{n-|C|-1}$, whenever $|C|$ grows by at most 2. Hence,
when the size of $C$ reaches $n - l \varepsilon^{-1} - 2$, we have
\begin{eqnarray*}
\dd &\geq& \sum_{i=2}^{n/2-\frac{1}{2\varepsilon}l-1}
\frac{1-2\varepsilon}{n-2i-1}\\
&\geq& \frac{1-2\varepsilon}{2}\sum_{i=4}^{n - l \varepsilon^{-1} -
2} \frac{1}{n-i-1}\\
&\geq& \frac{1-2\varepsilon}{2} \left(\sum_{i=1}^{n-5} \frac{1}{i} -
\sum_{i=1}^{l \varepsilon^{-1}} \frac{1}{i} \right)\\
&\geq& \frac{1-3\varepsilon}{2}\left(\log n - \log (l
\varepsilon^{-1}) \right)\\
&\geq& l,
\end{eqnarray*}
which concludes the proof of the theorem. {\hfill $\Box$}

\medskip

\section{Concluding remarks and open problems}
\label{sec::openprobs}

Recently, the approach we used to prove Theorem~\ref{planarity} was enhanced~\cite{AJSS}, and the error term was improved to a constant.

It was proved in Theorem~\ref{coloring} that $\tau_E({\cal NC}_n^2)
= \frac{n^2}{8} + \Theta (n)$. For $k \geq 3$, we know only the
simple bounds $\frac{(k-1)n^2}{4k} \leq \tau_E({\cal NC}_n^k) \leq
\frac{1}{2} {n \choose 2}$. Here the lower bound follows from
Tur\'an's Theorem and Observation~\ref{obs} and the upper bound is
the consequence of Enforcer being able to win. It would be
interesting to close, or at least reduce, the gap between these
bounds. It seems reasonable that, as in the case $k=2$, the truth is
closer to the lower bound, and maybe $\tau_E({\cal NC}_n^k) =
(1+o(1))\frac{(k-1)n^2}{4k}$, for every $k \geq 3$.

% In Theorem~\ref{posdeg=connect} it was shown that $\tau_E({\cal T}_n)$
% and $\tau_E({\cal D}_n)$ are ``quite close to each other''. This is
% reminiscent of the well-known property of random graphs, that the
% hitting time of being connected and the hitting time of having
% minimum positive degree are a.s.\ the same, and it motivates us to
% raise the following conjecture.
% \begin{conjecture} \label{hittingtime}
% $\tau_E({\cal D}_n) = \tau_E({\cal T}_n)$.
% \end{conjecture}

In Question~\ref{hittingtime} we ask whether $\tau_E({\cal D}_n) =
\tau_E({\cal T}_n)$ holds for sufficiently large $n$. It would be
interesting to consider related families, with a similar random
graph hitting time, like the hypergraph ${\mathcal M}_n$ of perfect
matchings or that of Hamilton cycles ${\mathcal H}_n$, and obtain
estimates on their minimum Avoider-Enforcer move number
$\tau_E({\cal M}_n)$ and $\tau_E({\cal H}_n)$.

\end{document}